\documentclass{article}

\usepackage{arxiv}

\usepackage[utf8]{inputenc} 
\usepackage[T1]{fontenc}    
\usepackage{hyperref}       
\usepackage{url}            
\usepackage{booktabs}       
\usepackage{amsfonts}       
\usepackage{nicefrac}       
\usepackage{microtype}      
\usepackage{lipsum}
\usepackage{amsmath}
\usepackage{optidef}
\usepackage{graphicx}
\usepackage{fancyvrb}
\usepackage{xcolor}
\graphicspath{ {./images/} }

\title{Identification of Time Delays in COVID-19 Data}

\author{
 Nicola Guglielmi \\
  Division of Mathematics\\
  Gran Sasso Science Institute\\
  Viale Francesco Crispi 7 \\
  L'Aquila, AQ 67100, Italy \\
  \texttt{nicola.guglielmi@gssi.it} \\
   \And
 Elisa Iacomini \\
  Institut fur Geometrie und Praktische Mathematik\\
  RWTH Aachen University\\
  Templergraben 55\\
  Aachen, 52062, Germany \\
  \texttt{iacomini@igpm.rwth-aachen.de} \\
  \And
 Alex Viguerie \\
  Division of Mathematics\\
  Gran Sasso Science Institute\\
  Viale Francesco Crispi 7 \\
  L'Aquila, AQ 67100, Italy \\
  \texttt{alexander.viguerie@gssi.it} \\
}

\begin{document}
\maketitle
\begin{abstract}
COVID-19 data released by public health authorities features the presence of notable time-delays, corresponding to the difference between actual time of infection and identification of infection. These delays have several causes, including the natural incubation period of the virus, availability and speed of testing facilities, population demographics, and testing center capacity, among others. Such delays have important ramifications for both the mathematical modeling of COVID-19 contagion and the design and evaluation of intervention strategies. In the present work, we introduce a novel optimization technique for the identification of time delays in COVID-19 data, making use of a delay-differential equation model. The proposed method is general in nature and may be applied not only to COVID-19, but for generic dynamical systems in which time delays may be present.
\end{abstract}


\section{Introduction}
The outbreak of COVID-19 in 2020 and into 2021 has led to a surge in interest in the mathematical modeling of the COVID-19 epidemic, as well as the mathematical modeling of epidemics generally. These models have taken many types. A full review of the relevant literature is beyond the current work. However, some general classes of models include data-driven and machine learning approaches \cite{linka2020safe, bhouri2021covid, jha2020bayesian, barros2021dynamic, viguerie2021coupled, wang2020system}, models based on partial differential equation (PDE) systems \cite{bertaglia2021hyperbolic, bertaglia2021hyperbolic2, albi2021kinetic, viguerie2021simulating, viguerie2020diffusion, grave2021adaptive, grave2021assessing, bertrand2021least}, agent-based models \cite{zohdi2020agent}, and models based on ordinary differential equation (ODE) systems. This last category is by far the most common such model, with such articles numbering in the thousands. Some works of this type include e.g. \cite{gatto2020spread, Ferguson2020, Gatto202004978, remuzzi2020covid, calafiore2020modified, loli2020monitoring, choi2020estimating, ivorra2020mathematical, parolini2021suihter}. We note that the categories listed above are not necessarily strict delimiters; indeed, many of the above cited works, as well myriad others, have characteristics of multiple of different model classes. One commonality that nearly all the presented models share is that they are, in some way, compartmental models in the SIR family owing to the seminal work of Kermack and Mackendrick \cite{kermack1927contribution} (see also: \cite{iannelli2015introduction, MurrayI} ).
\par One of the most important aspects of COVID-19 modeling deals with the presence of \textit{time delays} in the data. Time delays certainly manifest themselves via the natural incubation period of the disease; this is reflected in the employed models commonly incorporating an \textit{exposed} compartment, to model the period after exposure to the virus but before symptomatic infection \cite{linka2020safe, bhouri2021covid, viguerie2020diffusion, gatto2020spread}. Such an approach is a reasonable one, though not without difficulties. In particular, the exposed compartment is not necessarily easily known from data, requiring its estimation. This may be done, for example, through rule-of-thumb estimates based on the infected population, which result in high uncertainty \cite{grave2021assessing, grave2021adaptive, albi2021kinetic}. Compounding this problem further are the well-known presence of asymptomatic patients, who may spread the disease while never showing appreciable symptoms, and, at least in the earlier stages of the pandemic, were almost certainly under-counted \cite{albi2021kinetic, viguerie2020diffusion}.
\par In lieu of an exposed compartment, in \cite{guglielmi2021delay, viguerie2021coupled}, the authors employed a delay-differential equation (DDE) model, in which the incubation period was modeled with a time-delay term, rather than a separate compartment. Aside from providing the advantage of eliminating the potentially troublesome exposed compartment, such a modeling approach eliminates the assumption that the sojourn time of the incubation period is exponentially distributed, which is well-known to be a questionable assumption in some settings \cite{iannelli2015introduction}. Other works incorporating DDE models into the analysis of COVID-19 include \cite{dell2020solvable, devipriya2021seir, kumar2020analysis}, while references for more general epidemiological models of this type may be found in \cite{iannelli2015introduction, MurrayI, takeuchi2000global, brauer2012mathematical, forde2005delay, buonomo2008global}.
\par Aside from the presence of incubation periods, delays in data may result from other sources: additional time from initial onset of symptoms to when a person decides to be tested, the availability of testing, the speed of test result processing, testing center capacity, and other factors may have a large influence on when a given case is properly identified in the data \cite{sarnaglia2021correcting}. The presence of such lags is important, particularly when attempting to model measured data, as such processes are intrinsic to the data \cite{bastos2019modelling}. Hence, when attempting to properly evaluate different intervention strategies, for instance, such delays must be considered in the modeling process. A model designed to fit measured parameters, while not taking into account the delays in measurement, may fail provide a wholly accurate description of the dynamics. In this sense, we may view such delays as no longer time from exposure to \textit{infection}, but rather time from exposure to \textit{identification}.
\par In the present work, we seek to rigorously quantify this time delay effect in data. In order to accomplish this task, we first introduce, as a forward problem, a delayed-SIRD (susceptible-infected-recovered-deceased) model, incorporating a range of weighted possible time delays. We then employ an optimization process against measured data over the different weights, using the output of the optimization process to evaluate the presence of the different time delays present within the data. We note that the introduced method, while applied in the current work to COVID-19, is general in nature and such a technique may be employed in any situation in which a dynamical system and its associated measurements may exhibit time-delayed dynamics.
\par The article is outlined as follows. We first introduce the employed delay-differential equation model and discuss some of its characteristics. We then introduce and describe the optimization process for the identification of the relevant time delays. We then validate the proposed technique on an interval of data for the COVID-19 outbreak in Italy, demonstrating its effectiveness under several different modeling assumptions. Finally, we conclude with a summary of our main findings and possible directions for future research in this area.
\section{Delay differential equation model}
For the model used in the present work, we make the following assumptions: 
\begin{enumerate}
    \item That the time scales considered are sufficiently short such that we do not need to consider new births or non-COVID-19 deaths;
    \item That the time scales considered are sufficiently short such that we do not need to consider waning immunity;
    
    \item That probability of contact and transmission remains constant throughout the infectious period;
    \item That recovery and mortality rates follow an exponential distribution;
\end{enumerate}
We acknowledge that some of these assumptions are simplifications; however, for the scope of the current work, we feel they are nonetheless reasonable. 

Differently from \cite{guglielmi2021delay}  we model here the interaction with infected individuals by a convolution term
\begin{equation} \label{eq:convol}
\int\limits_{\sigma_{\min}}^{\sigma_{\max}} i(t-\sigma) w(\sigma) d \sigma 
\end{equation}
where $\sigma_{\min}$ and $\sigma_{\max}$ represent the extrema of the latency
period and $w(\sigma) \ge 0$ is a weight function (or a distribution) such that
\[
\int\limits_{\sigma_{\min}}^{\sigma_{\max}} w(\sigma) d \sigma = 1.
\]
We note that this term resembles the renewal function described first in \cite{kermack1927contribution}, from which the classical SIR model is then derived (see also \cite{breda2012formulation}). In the model discussed in \cite{guglielmi2021delay} we considered 
(with $\delta(\cdot)$ a Dirac-$\delta$ distribution)
\[
w(\sigma) = \delta(\sigma - \hat{\sigma}) 
\]
for a certain $\hat\sigma \in \left[ \sigma_{\min}, \sigma_{\max} \right]$ a representative constant delay identifying the latency period of the infection.
Here instead we aim to optimize the function $w$ in terms of the avalable measured data, allowing a certain variability of the latency, depending on the features of 
the population,  
In order to approach the problem numerically we discretize \eqref{eq:convol} by
\begin{equation*}
\int\limits_{\sigma_{\min}}^{\sigma_{\max}} i(t-\sigma) w(\sigma) d \sigma \approx
\sum_{j=1}^k w_j i(t-\sigma_j)
\end{equation*}
with $\{ \sigma_j \}$ a uniform partition of the interval 
$\left[ \sigma_{\min}, \sigma_{\max} \right]$.

In this way we can make use of an optimization procedure to estimate the weights,
and thus the memory effect, in terms of the observed data. As we will see the
discrete profile does not differ qualitatively from a bell-like function. This
indicates that there is a prevalent constant delay in the model and this justifies
the approximation considered in \cite{guglielmi2021delay}. The advantage is that in this way we are able to estimate numerically the most suitable value of a possibly
constant delay to be used in the model.

\medskip

Under the above assumptions, we are led to consider the following delayed-SIRD system of differential equations with delays $\sigma_j$, where $j=1:k$, with $k$ the number of considered delays:
\begin{align}
\label{eq1DelODE} \dot{s}(t) &=  - \beta s(t)\sum_{j=1}^k w_j i(t-\sigma_j),  \\
\label{eq2DelODE} \dot{i}(t) &= \beta s(t)\sum_{j=1}^k w_j i(t-\sigma_j)  - \left(\phi_d + \phi_r\right)\sum_{j=1}^k w_j i(t-\sigma_j), \\
\label{eq3DelODE} \dot{r}(t) &= \phi_r \sum_{j=1}^k w_j i(t-\sigma_j)    ,\\
\label{eq4DelODE} \dot{d}(t) &= \phi_d \sum_{j=1}^k w_j i(t-\sigma_j).
\end{align}
It was shown in \cite{guglielmi2021delay} that one may expect the model given by \eqref{eq1DelODE}-\eqref{eq4DelODE} to exhibit stable behavior provided that:
\begin{align}\label{stabilityCondition}
\phi_d + \phi_r < \frac{\pi}{2 \sigma},    
\end{align}
which in the present work is easily extended to: 
\begin{align}\label{stabilityCondition2}
\phi_d + \phi_r < \frac{\pi}{2 \sigma_j},    
\end{align}
for each $j$. 
\par We note additionally that above we have provided the \textit{density-dependent} version of the model for generality. One may easily obtain the \textit{frequency-dependent} version of \eqref{eq1DelODE}-\eqref{eq4DelODE} by instead considering $\widetilde{\beta} = \beta/N$, where $N$ is the total population, in lieu of $\beta$ in \eqref{eq1DelODE}, \eqref{eq2DelODE}. Such a change will not affect the stability bounds \eqref{stabilityCondition}-\eqref{stabilityCondition2}, though we note that the units of $\tilde{\beta}$ become Days$^{-1}$. 

\begin{table}
\begin{center}
\begin{tabular}{ |c|c|c| } 
\hline
Parameter/variable &  Name  &  Units   \\
\hline\hline

$s$ & Susceptible individuals & Persons \\ \hline
$i$ & Infected individuals & Persons \\ \hline
$r$ & Recovered individuals & Persons \\ \hline
$d$ & Deceased individuals & Persons \\ \hline
$\beta$ & Contact rate & Persons$^{-1} \cdot $ Days$^{-1}$  \\ \hline
$\phi_r$ & Recovery rate & Days$^{-1}$  \\ \hline
$\phi_d$ & Mortality rate & Days$^{-1}$  \\ \hline
$\sigma_j$, $j=1:k$ & Time-delay $j$ & Days   \\ \hline
$w_j$, $j=1:k$ & Time-delay weight $j$ & Dimensionless   \\ \hline
\end{tabular}
\caption{Relevant parameters, parameter names, and units for the equations.}
  \label{tab:1DParametertTable} 
\end{center}
\end{table}

\section{Optimization procedure and methods}
We seek to use \eqref{eq1DelODE}-\eqref{eq4DelODE} identify the time delays $\sigma_j$ present in COVID-19 data. To this end, we solve the inverse problem for the delay weights $w_j$. Denote the measured susceptible, infected, recovered, and deceased (from COVID-19) populations as $\widetilde{s},\,\widetilde{i},\,\widetilde{r},\,\widetilde{d}$, respectively. We then define the following optimization problem: For a given $\beta$, $\phi_r$, $\phi_d$, $\sigma_j$, $j=1:k$ and $\boldsymbol{w}=\lbrack w_1,\,w_2,\,...,\,w_k\rbrack^T$:
\begin{argmini!}|l|[3]
{\boldsymbol{w}}{(i-\widetilde{i})^2 + (d-\widetilde{d})^2}
{}{}
\addConstraint{\sum_{j=1}^k w_j = 1}
\addConstraint{0\leq w_j \leq 1\,\,\forall j,}
\addConstraint{s,\,i,\,r,\,d \text{ solve } \eqref{eq1DelODE}-\eqref{eq4DelODE}}{}.
\end{argmini!}
\par The novelty of the proposed method above lies in its treatment of the delay term as a weighted sum of different delays. By considering the system in this manner, it is possible to obtain an optimization problem for the time delay that is linear in the delay term. 
This is important, as attempting to optimize directly on the $\sigma_j$ gives a highly nonlinear optimization problem that is likely untractable. 
This method of discretizing the time lags and optimizing for the weights as a means of time-delay identification is, to the authors' knowledge, novel. We note that this approach is highly general and is not restricted to models of COVID-19 or epidemiological models in general. 
{Moreover, the accuracy and the detailed description of the delay depend on how $k$ is chosen. This allows this method to be very applicable for other models as well. } 
\par We note that the model \eqref{eq1DelODE}-\eqref{eq4DelODE} depends heavily not only on the time delays $\sigma_j$ and their respective weights $w_j$, but also on the contact rate $\beta$, recovery rate $\phi_r$, and mortality rate $\phi_d$. In order to obtain robust estimates for the $w_j$ through the optimization procedure defined above, we resort to a Monte Carlo method. 
\par We consider the parameters $\beta$, $\phi_r$, $\phi_d$ as independent random variables sampled from Gaussian distributions. The mean of the respective distributions is obtained via an empirical parameter fitting, with the standard deviation of each distribution taken to be 5\% the value of the mean. We then run the optimization procedure (8a)-(8d) a total of 1000 times, corresponding to 1000 different parameter samplings. We then identify frequency with which different $w_j$ appear in the optimization, as well as the mean value of each $w_j$ over the sampled trial.

\par We note as well, that a limitation of our implementation is its equal treatment of time delays across compartments. In general, we consider the infection and deceased data to be the most reliable such data and is emphasized in our calculations. {As the easiest quantity to measure, the deceased compartment is assumed the most accurate. Further, the as the rate of new infections is the driving point behind policy measures, it is considered the second most accurate measure for our purposes. }
However, one may modify (8a) to include other compartments, or different subsets of the compartment set. 

\par In terms of implementation, there are several important points to note. We define the optimization procedure using the MATLAB routine \verb|fmincon|. In order to solve the delay differential equation system \eqref{eq1DelODE}-\eqref{eq4DelODE}, we use the MATLAB function \verb|dde23|. We note that \verb|dde23| uses adaptive grid points, and hence the temporal points of the solution obtained using \verb|dde23| do not correspond to those in our dataset in general. To circumvent this, we use \verb|deval| to evaluate the \verb|dde23| solution at the desired temporal points. 

\section{Results}
For our measured data, we consider the COVID-19 outbreak in Italy, as reported by the newspaper Il Sole 24 \cite{Lab24}. We consider the dates ranging from September 11, 2020 to February 7, 2021. This date range was chosen to represent a 150-day span, in which different governmental policies were enacted, widespread testing was available, and before large-scale population vaccination. This choice was intentional, as we wish to avoid as many confounding factors as possible. We consider the following distributions for the parameters
\begin{align}
    \beta \sim N\left(.1131/n_0,\,.05(.1131/n_0) \right) \\
    \phi_d \sim N\left(1/940,\,.05(1/940)\right) \\
    \phi_r \sim N(1/24,\,.05(1/24)),
\end{align}
where $N(\mu,\sigma)$ denotes a normal distribution with mean $\mu$ and standard distribution $\sigma$ and $n_0$ denotes the initial living population. For $t>73$, where $t_0=0$, we replace $\beta$ by $\beta/3$ to model the introduced government restrictions, designed to curb the virus spread. The values of $\phi_d$ and $\phi_r$ remain unchanged throughout the considered interval. The mean values of these parameters were obtained via a preliminary parameter tuning. For the delays, we considered 12 possible values, from $\sigma=2$ days to $\sigma=35$ days, with 3 days of spacing between the different days. We also performed some experiments (not shown) using different discretizations over the same range (for instance, 10 possible values) to show the robustness of the optimization procedure, which did not appreciably change the conclusions drawn. 
\par To examine the possible differences in delays among the different compartments, we perform the optimization four times: once considering the $i$ and $d$ compartments together in the objective function, and for the $i$, $r$ and $d$ compartments, considering each compartment individually. For each ensemble, we report the minimum, maximum, and mean $L^2$ error over the simulation ensemble when compared to the measured data, as well as the mean value and frequency of the different weights $w_j$ corresponding to the delays $\sigma_j$.

\par
\begin{figure}

\includegraphics[width=\linewidth]{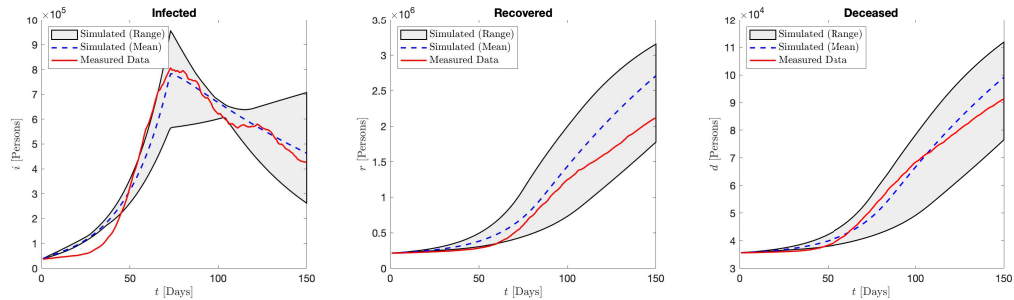}\caption{Aggregated results of the infected, recovered, and deceased compartments over the range of simulations when optimizing over the $i$ and $d$ compartments. We see good agreement across all compartments with the measured data, and particularly strong agreement (less than 10\% in $L^2$ norm) in the deceased compartment.}\label{fig:SimResultsID}
\end{figure}

 The aggregated results of the simulations in terms of the infected, recovered, and deceased compartments, when optimizing over the $i$ and $d$ compartments are shown in Fig. \ref{fig:SimResultsID}. We see good qualitative agreement with the measured and simulated data across the three compartments of interest, with the major relevant trends being captured. In table \ref{tab:ErrorTableID}, we report the mean, minimum, and maximum of each error value computed over the simulation ensemble. We see particularly strong agreement in the susceptible and deceased compartments, reasonable agreement in the infected compartment (especially considering its rapid-changing trajectory), and somewhat worse performance in the recovered compartment. Given that the $r$ compartment was not included in the optimization, this is unsurprising.

\begin{table}
\begin{center}
\begin{tabular}{ |c|c|c|c| } 
\hline
Quantity & Mean &  Min.  &  Max.    \\
\hline\hline
$\|s-\widetilde{s}\|_{L^2}/\|\widetilde{s}\|_{L^2}$ & .0051 & .0022 & .0104 \\ \hline
$\|i-\widetilde{i}\|_{L^2}/\|\widetilde{i}\|_{L^2}$ & .1101 & .0798 & .2757 \\ \hline
$\|r-\widetilde{r}\|_{L^2}/\|\widetilde{r}\|_{L^2}$ & .2602 & .071 & .5913 \\ \hline
$\|d-\widetilde{d}\|_{L^2}/\|\widetilde{d}\|_{L^2}$ & .0695 & .0249 & .2063 \\ \hline
\end{tabular}
\caption{Error table: relative errors in $L^2$ norm against measured data for the different compartments over 1000 simulations, optimizing over $i$ and $d$.}
  \label{tab:ErrorTableID} 
\end{center}
\end{table}

\begin{table}
\begin{center}
\begin{tabular}{ |c|c|c|c| } 
\hline
Quantity & Mean &  Min.  &  Max.    \\
\hline\hline
$\|s-\widetilde{s}\|_{L^2}/\|\widetilde{s}\|_{L^2}$ & .0051 & .0022 & .0108 \\ \hline
$\|i-\widetilde{i}\|_{L^2}/\|\widetilde{i}\|_{L^2}$ & .1116 & .0807 & .2681 \\ \hline
$\|r-\widetilde{r}\|_{L^2}/\|\widetilde{r}\|_{L^2}$ & .2561 & .0650 & .6054 \\ \hline
$\|d-\widetilde{d}\|_{L^2}/\|\widetilde{d}\|_{L^2}$ & .0684 & .0266 & .1659 \\ \hline
\end{tabular}
\caption{Error table: relative errors in $L^2$ norm against measured data for the different compartments over 1000 simulations, optimizing only over $i$.}
  \label{tab:ErrorTableIOnly} 
\end{center}
\end{table}

\begin{table}
\begin{center}
\begin{tabular}{ |c|c|c|c| } 
\hline
Quantity & Mean &  Min.  &  Max.    \\
\hline\hline
$\|s-\widetilde{s}\|_{L^2}/\|\widetilde{s}\|_{L^2}$ & .0025 & .0016 & .0041 \\ \hline
$\|i-\widetilde{i}\|_{L^2}/\|\widetilde{i}\|_{L^2}$ & .2347 & .1580 & .4672 \\ \hline
$\|r-\widetilde{r}\|_{L^2}/\|\widetilde{r}\|_{L^2}$ & .0586 & .0394 & .1091 \\ \hline
$\|d-\widetilde{d}\|_{L^2}/\|\widetilde{d}\|_{L^2}$ & .0926 & .0362 & .1672 \\ \hline
\end{tabular}
\caption{Error table: relative errors in $L^2$ norm against measured data for the different compartments over 1000 simulations, optimizing only over $r$.}
  \label{tab:ErrorTableROnly} 
\end{center}
\end{table}

\begin{table}
\begin{center}
\begin{tabular}{ |c|c|c|c| } 
\hline
Quantity & Mean &  Min.  &  Max.    \\
\hline\hline
$\|s-\widetilde{s}\|_{L^2}/\|\widetilde{s}\|_{L^2}$ & .0035 & .0009 & .0079 \\ \hline
$\|i-\widetilde{i}\|_{L^2}/\|\widetilde{i}\|_{L^2}$ & .1464 & .0826 & .4992 \\ \hline
$\|r-\widetilde{r}\|_{L^2}/\|\widetilde{r}\|_{L^2}$ & .1903 & .0440 & .4769 \\ \hline
$\|d-\widetilde{d}\|_{L^2}/\|\widetilde{d}\|_{L^2}$ & .0430 & .0170 & .0864 \\ \hline
\end{tabular}
\caption{Error table: relative errors in $L^2$ norm against measured data for the different compartments over 1000 simulations, optimizing only over $d$.}
  \label{tab:ErrorTableDOnly} 
\end{center}
\end{table}

\par
\begin{figure}[!htb]
\minipage{0.5\textwidth}
  \includegraphics[width=\linewidth]{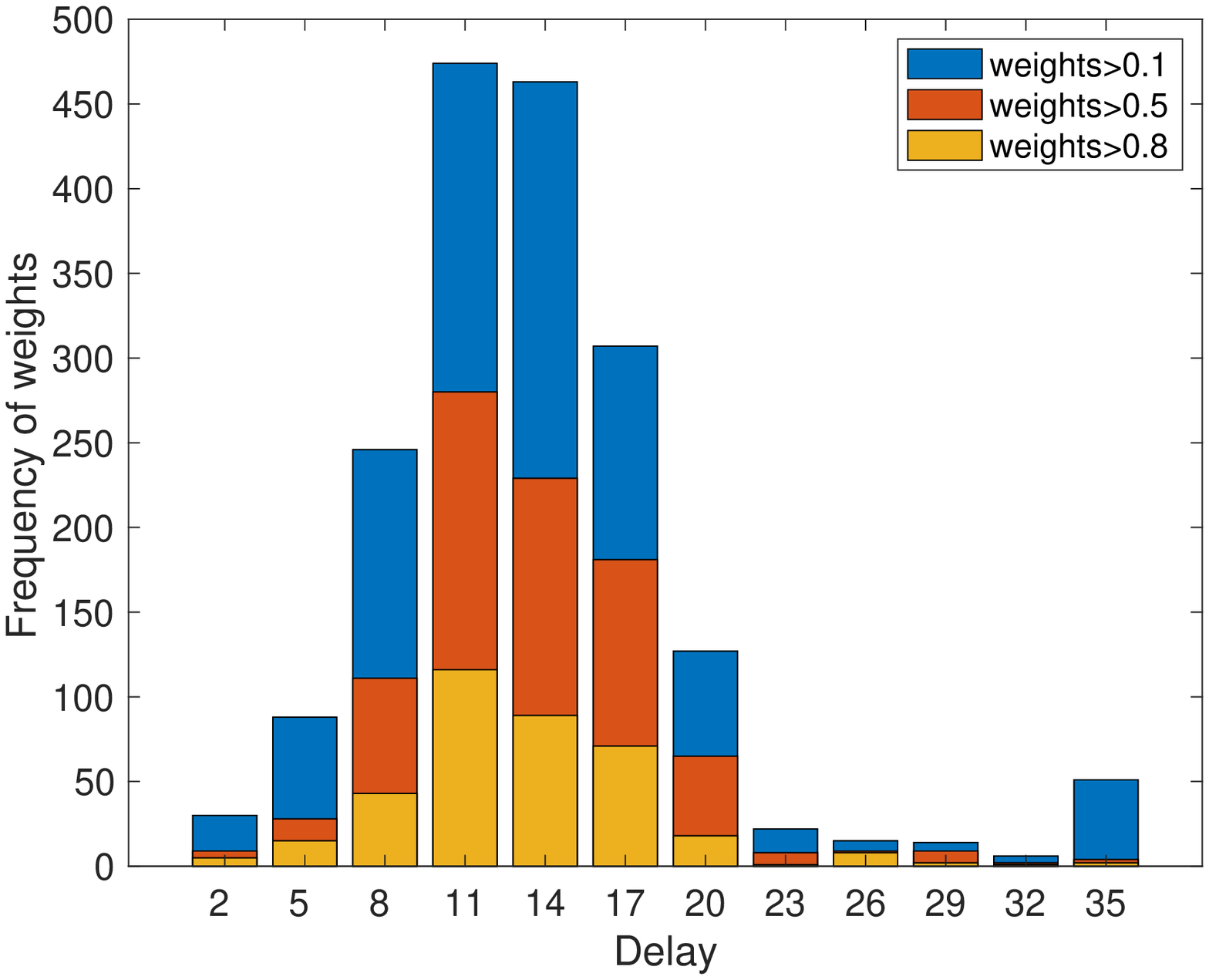}
 \endminipage\hfill
\minipage{.5\textwidth}%
  \includegraphics[width=\linewidth]{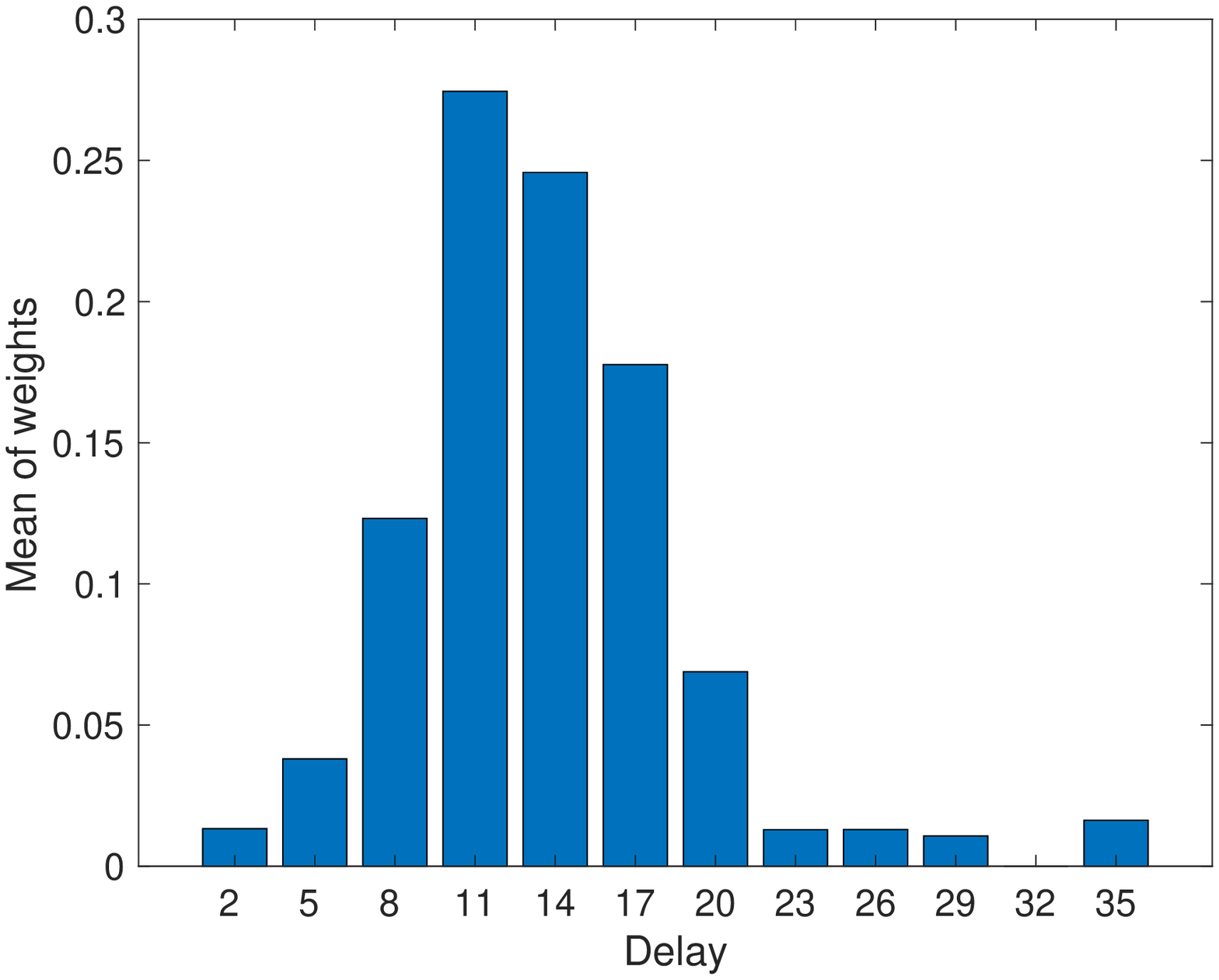}
\endminipage\caption{Left: the different values of the $w_j$ by frequency of occurrence when optimizing over $i$ and $d$. We see a strong concentration of weights in the range of $\sigma=$8 to $\sigma=$17 days. Right: the mean values of the weights over the simulation ensemble. These results mirror those of the frequency plot, showing a concentration between 8 and 17 days. For both frequency and mean, $\sigma=11$ appears as the most represented delay. } 
\label{fig:WeightsID}
\end{figure}

 Turning our attention to the influence of time delays, we report the mean and frequency of the different weights $w_j$ in Figure \ref{fig:WeightsID}. We see a heavy concentration of weights values, for both frequency and mean, in the range of $\sigma=8$ to $\sigma=17$ days. This is consistent with what has been reported in general media, as the incubation period of COVID is believed to have an average of 3-11 days \cite{lauer2020incubation}. When combined with natural delays coming from the testing and data reporting, a total time lag in the range of 8-17 days would appear consistent with this information.

\par
\begin{figure}

\includegraphics[width=\linewidth]{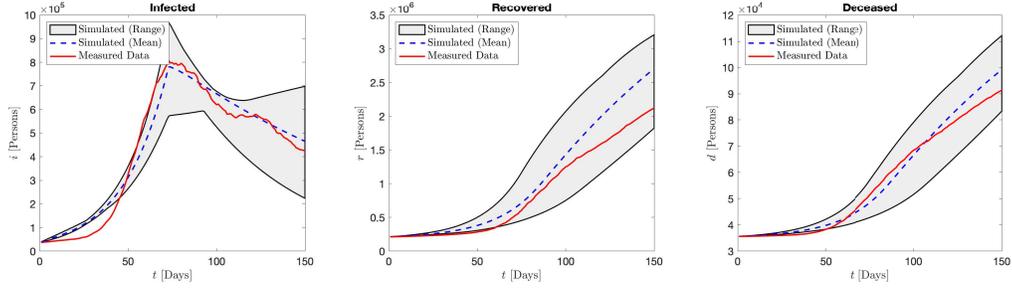}\caption{Aggregated results of the infected, recovered, and deceased compartments over the range of simulations when optimizing over the $i$ compartment. We see good agreement across all compartments with the measured data, and particularly strong agreement (less than 10\% in $L^2$ norm) in the deceased compartment.}\label{fig:SimResultsI}
\end{figure}

\begin{figure}

\includegraphics[width=\linewidth]{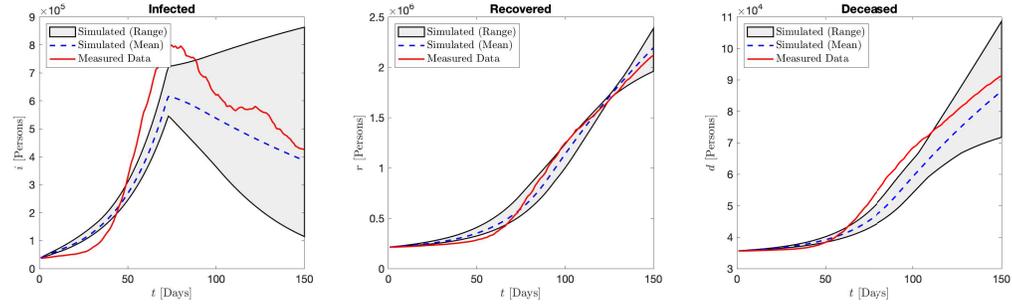}\caption{Aggregated results of the infected, recovered, and deceased compartments over the range of simulations when optimizing over the $r$ compartment. We see good agreement in the $r$ compartment; the agreement with the $d$ compartment remains acceptable, though somewhat lesser than other optimization choices. The error in the infected compartment is notably higher. This different error behavior is reflected in the optimization results; when optimizing over $r$, the optimal time delays are different than when optimizing over $i$ and $d$.}\label{fig:SimResultsR}
\end{figure}

\begin{figure}

\includegraphics[width=\linewidth]{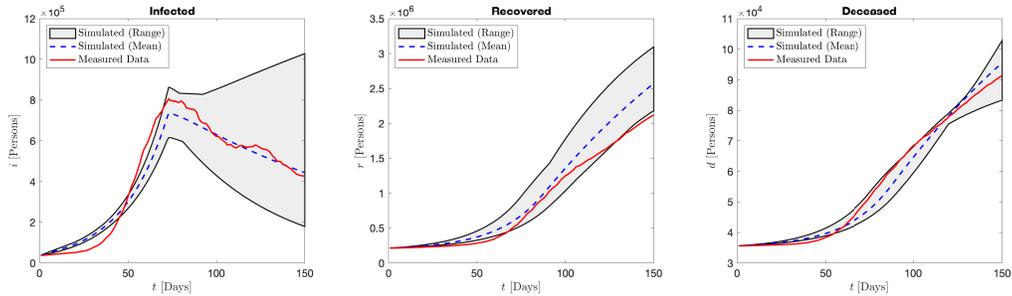}\caption{Aggregated results of the infected, recovered, and deceased compartments over the range of simulations when optimizing over the $d$ compartment. We see good agreement in the $d$ and $r$ compartments, with somewhat less close agreement with the $i$ compartment. }\label{fig:SimResultsD}
\end{figure}

 We then repeat this experiment three more times: once each for the $i$, $r$, and $d$ compartments considered individually in the optimization process. The error tables for the $i$, $r$ and $d$ compartments can be found in Table \ref{tab:ErrorTableIOnly}, \ref{tab:ErrorTableROnly}, and \ref{tab:ErrorTableDOnly} respectively. The corresponding simulation outputs for each case are plotted in Figures \ref{fig:SimResultsI}, \ref{fig:SimResultsR} and \ref{fig:SimResultsD}. The behavior is consistent with what one may intuitively expect; in each instance, the error behavior is minimized for the considered compartment, with the simulation variance greatly reduced for the isolated compartment. We note that one can observe from the relevant tables and figures that there is a much higher degree similarity between the $i$ and $d$ error behavior when compared to the $r$ compartment, which appears to act more independently.

\begin{figure}[!htb]
\minipage{0.5\textwidth}
  \includegraphics[width=\linewidth]{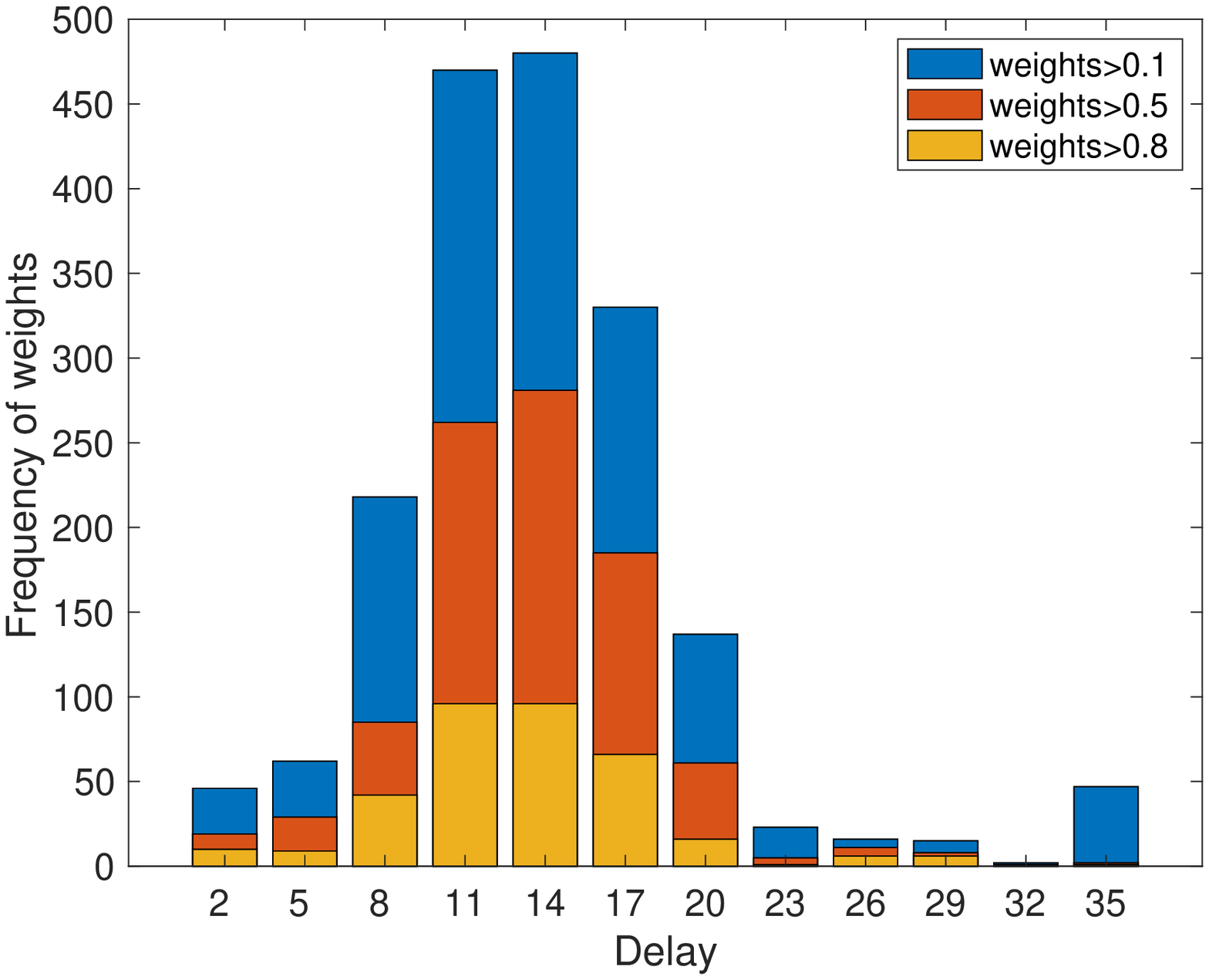}
 \endminipage\hfill
\minipage{.5\textwidth}%
  \includegraphics[width=\linewidth]{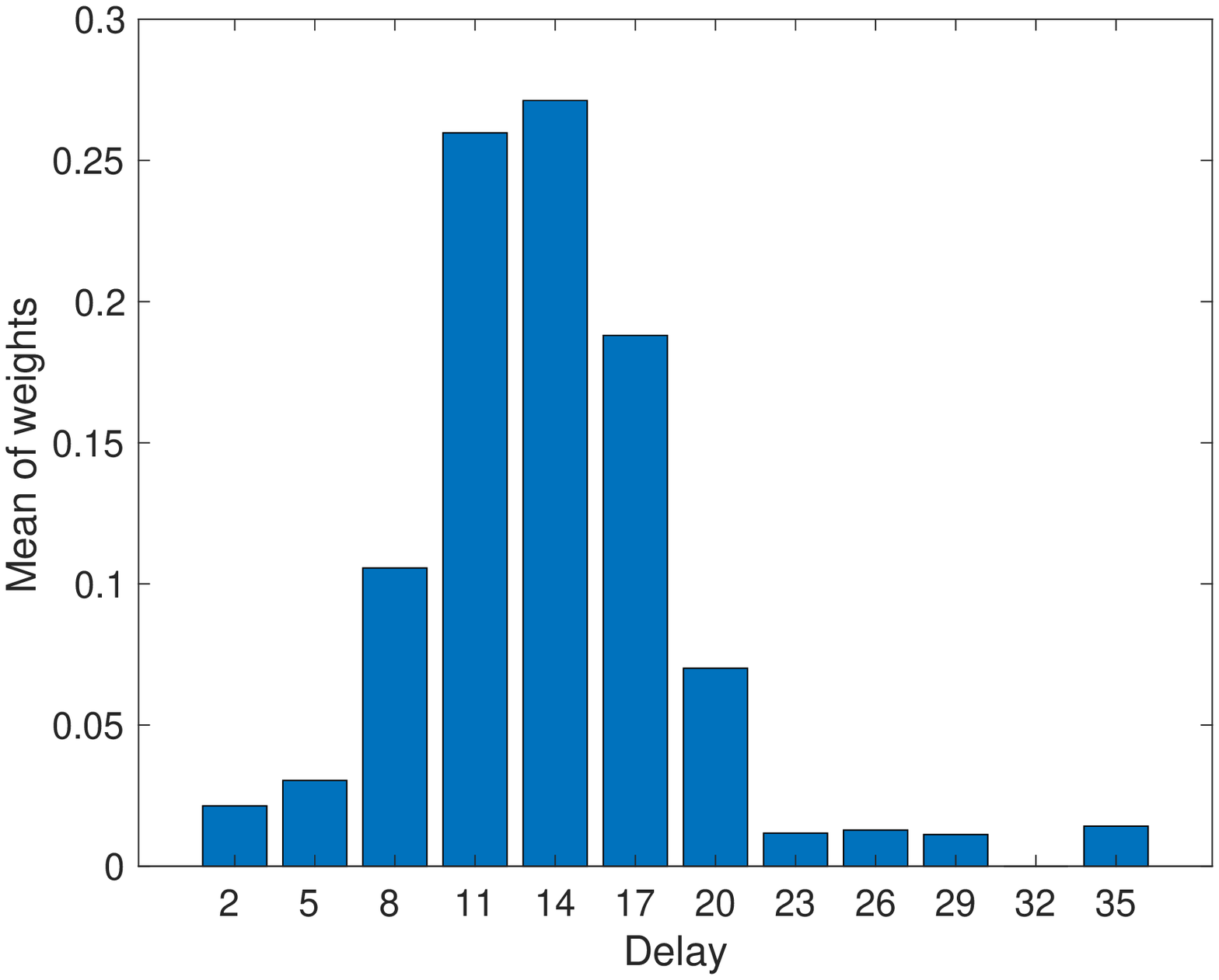}
\endminipage\caption{Left: the different values of the $w_j$ by frequency of occurrence when optimizing over $i$. We see a strong concentration of weights in the range of $\sigma=$8 to $\sigma=$17 days. Right: the mean values of the weights over the simulation ensemble. These results mirror those of the frequency plot, showing a concentration between 8 and 17 days. For both frequency and mean, $\sigma=14$ appears as the most represented delay.  } 
\label{fig:WeightsI}
\end{figure}

\begin{figure}[!htb]
\minipage{0.5\textwidth}
  \includegraphics[width=\linewidth]{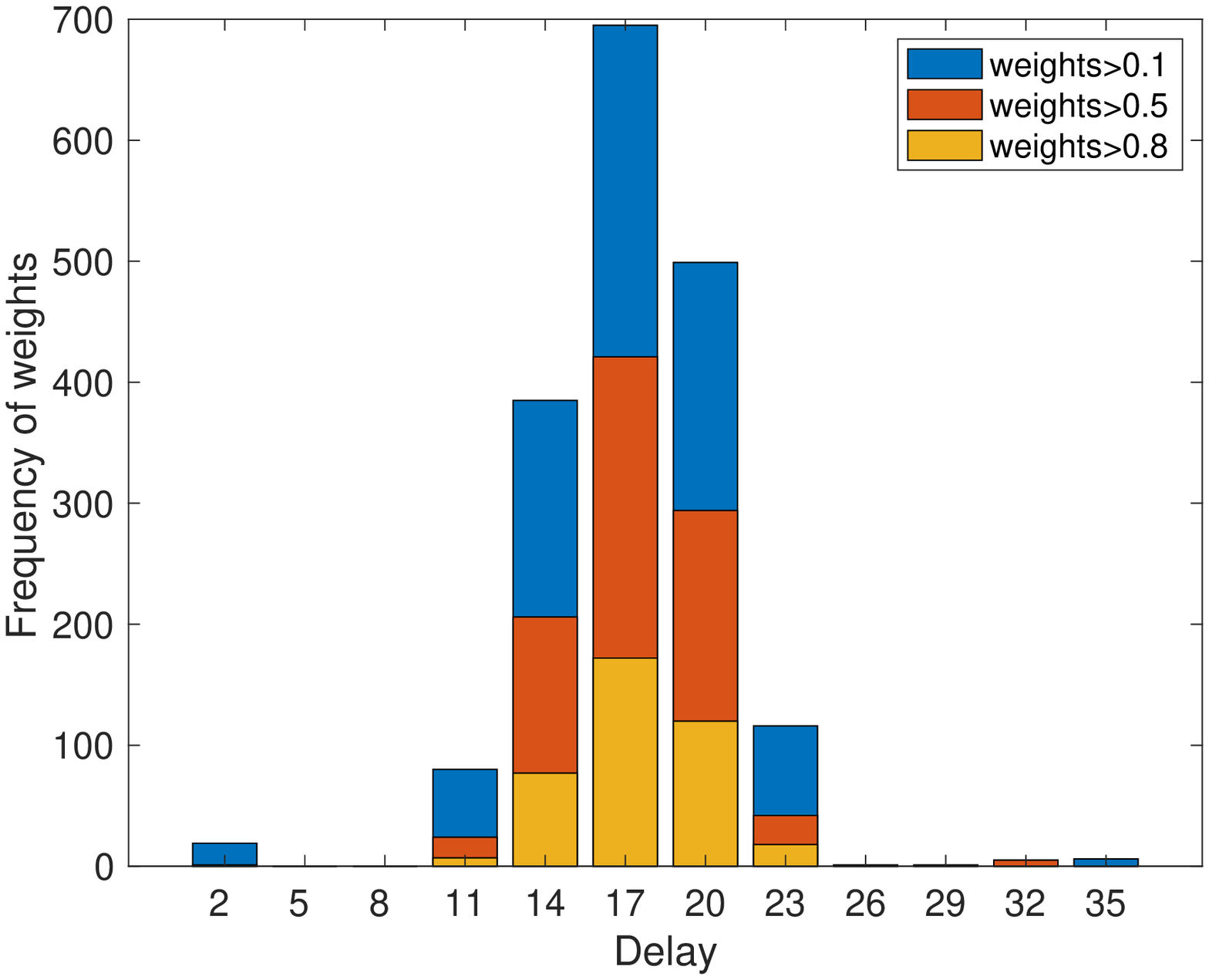}
 \endminipage\hfill
\minipage{.5\textwidth}%
  \includegraphics[width=\linewidth]{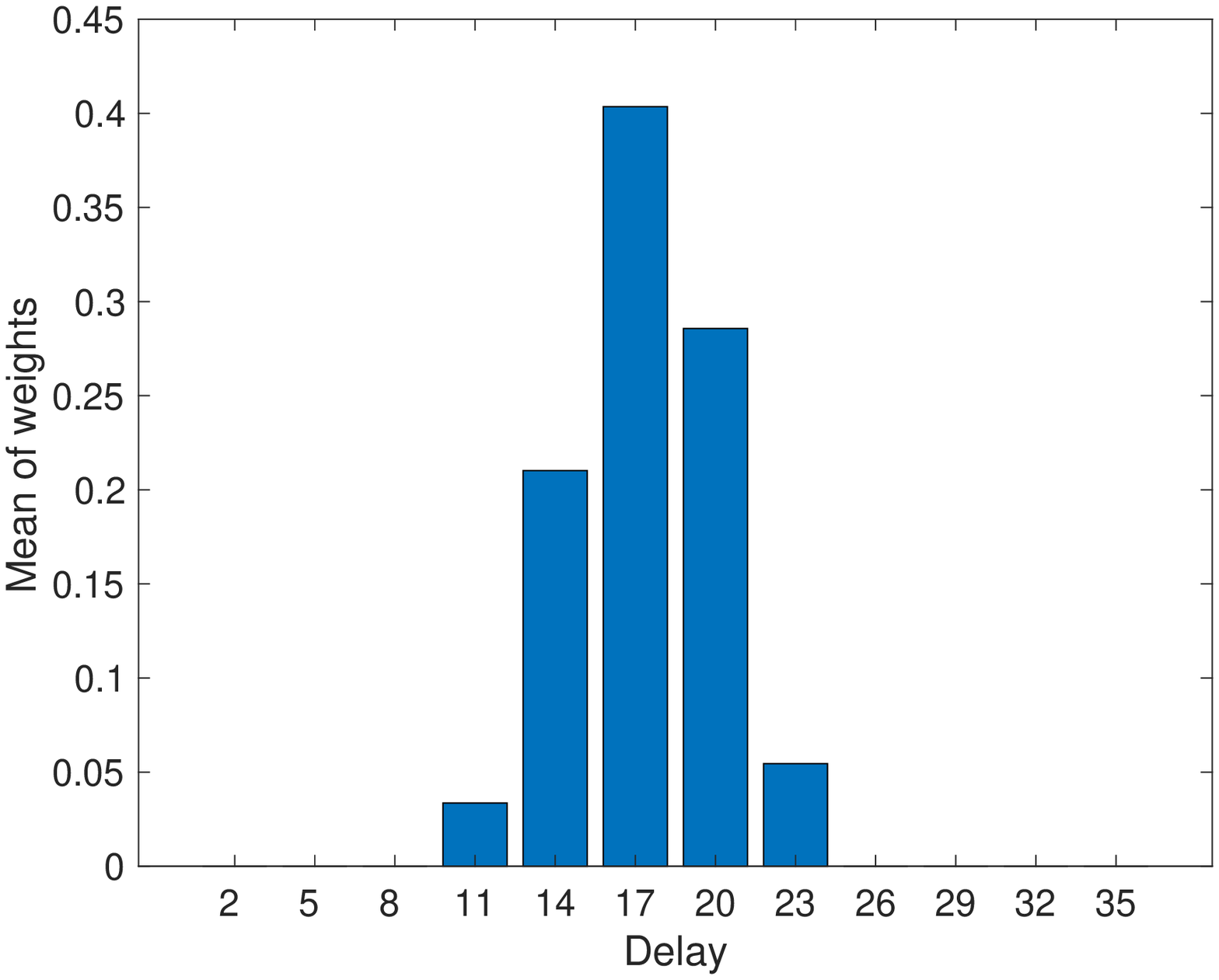}
\endminipage\caption{Left: the different values of the $w_j$ by frequency of occurrence when optimizing over $r$. We see a strong concentration of weights in the range of $\sigma=$14 to $\sigma=$20 days. Right: the mean values of the weights over the simulation ensemble. These results mirror those of the frequency plot, showing a concentration between 14 and 20 days. For both frequency and mean, $\sigma=17$ appears as the most represented delay. These results are less similar to those obtained when one optimizes over $i$ and $d$, and suggests different time-delays present in the $r$ compartment as compared to other compartments } 
\label{fig:WeightsR}
\end{figure}

\begin{figure}[!htb]
\minipage{0.5\textwidth}
  \includegraphics[width=\linewidth]{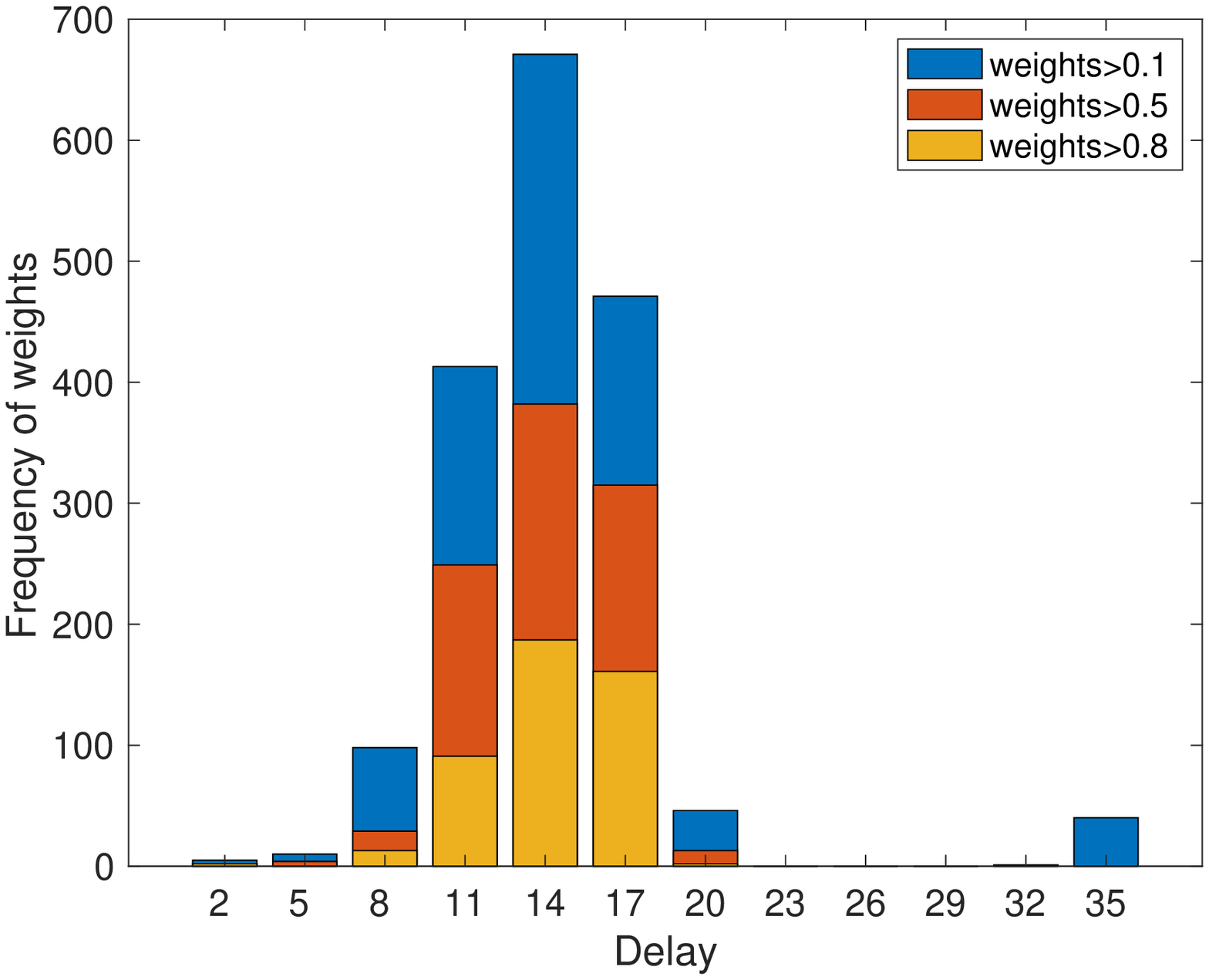}
 \endminipage\hfill
\minipage{.5\textwidth}%
  \includegraphics[width=\linewidth]{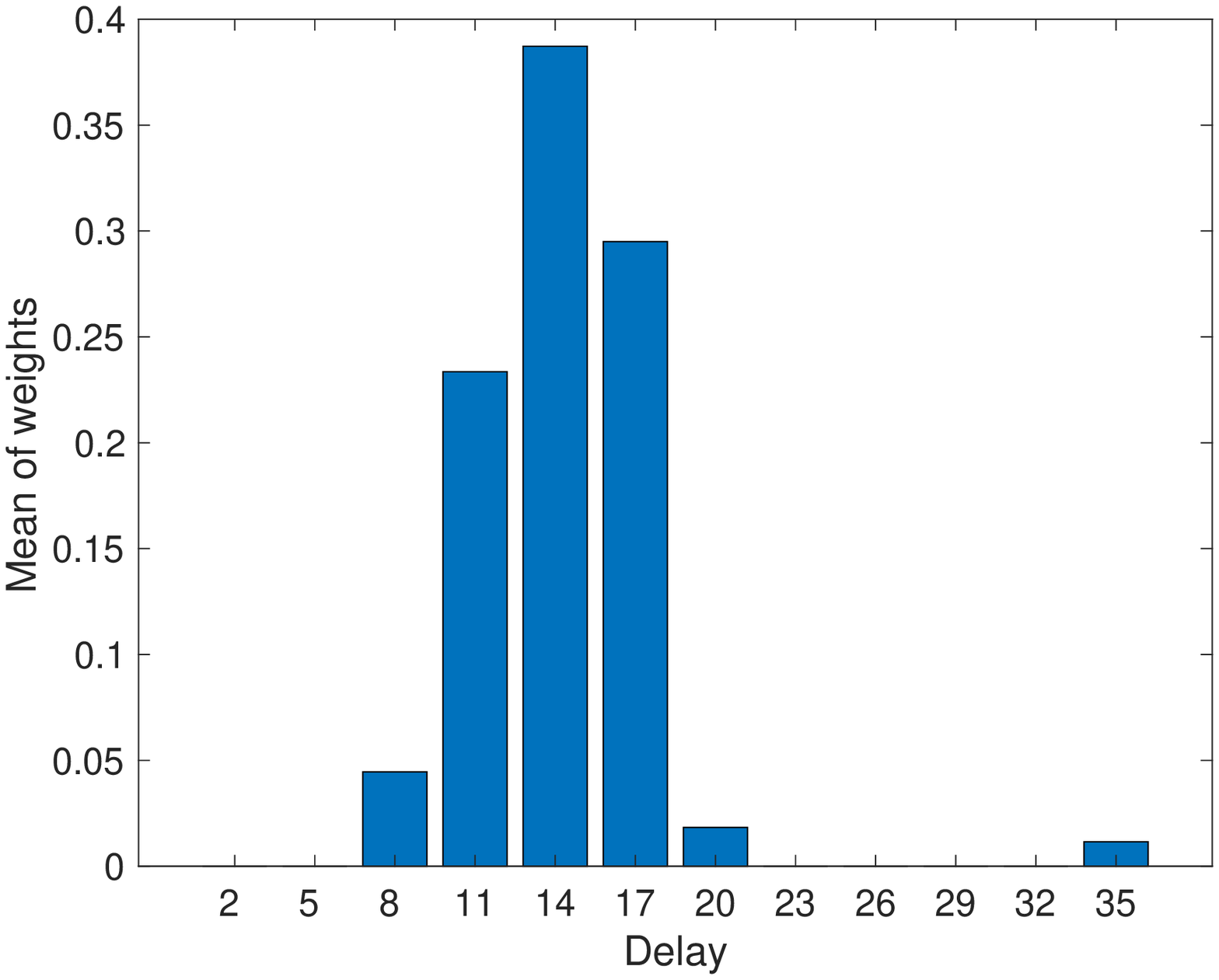}
\endminipage\caption{Left: the different values of the $w_j$ by frequency of occurrence when optimizing over $d$. We see a tight concentration of $\sigma$ in the range of 11 to 17 days. Right: the mean values of the weights over the simulation ensemble. These results mirror those of the frequency plot, showing a concentration between 11 and 16 days. For both frequency and mean, $\sigma=14$ appears as the most represented delay.} 
\label{fig:WeightsD}
\end{figure}
\par This is confirmed when looking at the frequency and mean distribution of the weights $w_j$ in Figures \ref{fig:WeightsI}, \ref{fig:WeightsR} and \ref{fig:WeightsD} for the $i$, $r$, and $d$ compartments respectively. In the case of the $i$ compartment, we see concentration between $\sigma$=8 and $\sigma$=17 days, similar to when $i$ and $d$ are considered jointly; however, for $i$ alone, $\sigma=$14 appears as the most dominant lag term, rather than $\sigma=$11. The behavior appears similar but even more pronounced when considering $d$ alone. In this instance $\sigma$=14 is again the dominant delay term; however, the grouping is tightly clustered between $\sigma=$11 and $\sigma=$17, with terms outside this range showing very little influence.
\par In contrast, when examining the case of isolated $r$ in \ref{fig:WeightsR}, we see a different distribution of delay terms. Indeed, in this instance, the dominant delay terms range from $\sigma=14$ to $\sigma=20$, and are very tightly concentrated within this range. This suggests a distinct and longer time scale when compared to the $i$ and $d$, and may explain why the error behavior for the $r$ compartment shows a distinctly different behavior in the plots and tables shown. 
\par In all, we can conclude three key takeaway points from these results:
\begin{enumerate}
    \item The time delays present in COVID-19 data range from $\sigma=8$ to $\sigma=20$ days, with the most likely value for $\sigma_j$ the $i$ and $d$ compartments in the range of $\sigma$=11 to $
    \sigma$=14 days.
\item The scale of the time delays across the compartments is not uniform. It is similar for new infections and for deceased persons, but occurs over a somewhat longer time scale for recovered individuals.{This is not surprising due to the different recovery rate and the inherent difficulty of measuring this compartment, given delays in testing and in the reporting of results.}
\item The time scales revealed by the optimization procedure shown here are longer than most estimates of the incubation period, suggesting that other delays are present in the data, and a period of 14 to 17 days should be considered when analyzing the effects of interventions. 
\end{enumerate}

\section{Conclusions and Future Directions}
\par In the present work, we have attempted to identify time delays present in COVID-19 data by means of newly-introduced method. This method consists of discretizing $k$ possible time delays  $\sigma_j$, $j=1:k$ and expressing these terms as a weighted sum $\sum_{j=1}^k w_j i(t-\sigma_j)$. This formulation then yields a an optimization problem which is linear in the weights $w_j$. After a preliminary parameter fitting, we then solve an ensemble of optimization problems, finding the weights $w_j$ which minimize the difference between simulated and measured data over different samplings of parameter values. We recover a distribution of delay terms for the different variables, allowing us to recover the inherent time-delays within the measured data.
\par We tested this approach on COVID-19 data in Italy, and found that the time-delays differ among different data points, but generally speaking range between 8 and 17 days for $i$ and $d$, with the most likely value between 11 and 14 days, and over a somewhat longer time scale (14-20 days) for the recovered compartment. This is consistent with our expectations given the incubation period of the disease and the natural time delays one may expect from the testing and data collection process. Our results confirm that proper evaluation of a given intervention for COVID-19 should allow for a sufficiently long time window to pass before conclusions are made.
\par There are several important directions for future work. Given the increased importance of waning vaccine efficacy over time, a similar analysis given vaccine data and data on breakthrough infections could be helpful in identifying similar such time-delay trends for this problem. While we have elected to avoid vaccines in the current work, preferring to establish a proof of concept on a more simplified problem, the approach offered is extremely general and may be applied to these important problems. Further, given the generality of the proposed optimization procedure, it may have wide applications outside of COVID-19 and epidemiology generally, and may be applied to any model where time delays are present, including applications in traffic modeling, pedestrian flow, remote sensing, and other fields. From the theoretical and methodological standpoint, an extension of the framework shown here for problems in which the delays are non-constant or state-dependent is also an important direction of future research.
\bibliographystyle{unsrt}  

\bibliography{references.bib}

\end{document}